\documentclass[a4paper]{amsart}
\usepackage[utf8x,utf8]{inputenc}
\usepackage{fancyhdr} 
\usepackage{amsfonts,amsmath,amssymb,amsthm}
\usepackage{amscd}
\usepackage{mathrsfs}
\usepackage{ucs,bbm} 
\usepackage{colonequals}
\usepackage[all]{xy}
\usepackage{verbatim}
\usepackage[backref=page]{hyperref}

\hypersetup{
     colorlinks   = true,
     citecolor    = red,
     linkcolor = blue,
     urlcolor = purple
}

\usepackage{enumerate}
\usepackage{xcolor}

\mathchardef\dash="2D

\newcommand{\dR}{\mathrm{dR}}

\newcommand{\calO}{\mathcal{O}}

\newcommand{\BB}{\mathbb{B}}
\newcommand{\CC}{\mathbb{C}}

\newcommand{\ZZ}{\mathbb{Z}}

\newcommand{\Fil}{\mathrm{Fil}}
\newcommand{\NN}{\mathbb{N}}
\newcommand{\GG}{\mathbb{G}}

\DeclareMathOperator{\Lie}{Lie}

\DeclareMathOperator{\Coker}{Coker}

\DeclareMathOperator{\dlog}{dlog}

\newcommand{\cvec}[1]{ \begin{pmatrix} #1 \end{pmatrix} }

\theoremstyle{plain}
\newtheorem{theorem}{Theorem}[section]

\theoremstyle{definition}

\theoremstyle{remark}
\newtheorem{remark} [theorem]{Remark}

\title{On log differentials of local fields}

\author{Nicola Mazzari}
\address{Università degli Studi di Padova}
\email{mazzari@math.unipd.it}

\begin{document}


\maketitle
\begin{abstract}
  We prove a  logarithmic version of  Fontaine's classic result on differentials of $\calO_{\bar K}$ over $\calO_K$. 
\end{abstract}



\section{Introduction}
Let $K$ be a $p$-adic field,  $\bar K$  an algebraic closure and $\CC_p=\hat{\bar{K}}$ its completion. We denote by $\calO_{K}$ and $\calO_{\bar K}$ their rings of integers. The module of K\"ahler differentials $\Omega=\Omega_{\calO_{\bar K}/\calO_{ K}}$ is well-known.
In particular  Fontaine \cite[Corollaire 1]{Fon:81a} gives the following identification
\[
	\CC_p(1)\cong V_p(\Omega)
\]
 where  $V_p(\Omega)$ is the rational Tate module of $\Omega$. The latter isomorphism is induced by the  map
 \[
 	\dlog_F:\bar K(1)=\bar K \otimes_{\ZZ_p} \ZZ_p(1)\rightarrow \Omega\ , \  \dlog_F(p^{-n}\otimes \underline{\epsilon})=\frac{d\epsilon_n}{\epsilon_n}\ ,
 \]
 where $\underline{\epsilon}=(\epsilon_n)_{n\ge 0}$ is a compatible sequence of $p$-powers roots of $1$, with $\epsilon_1\ne 1$.
 
 In this article we slightly modify the source of the previous map in order to compute a logarithmic version of $\Omega$ and its Tate module.
 
Fix $q\in p\calO_K\setminus \{0\}$. We denote by $E_q$ the Tate elliptic curve of parameter $q$: its group of  $\bar{K}$-points is  isomorphic to ${\bar K}^*/\langle q\rangle$ and its Tate module, denoted $\ZZ_p(1)^{\log}$, is an extension of the trivial Galois module $\ZZ_p$ by $\ZZ_p(1)$. Its divisible group is denoted by $\mu_{p^\infty}^{\log}$. 
Fix a sequence $\underline{q}=(q_n)_n$ of compatible $p$-power roots of $q$.  Then, $\ZZ_p(1)^{\log}=\ZZ_p\tilde{\epsilon}\oplus\ZZ_p\tilde{q}$ as a $\ZZ_p$-module, where $\tilde{\epsilon},\tilde{q}$ are the logarithms of $\underline{\epsilon},\underline{q}$. The Galois action is given by $\cvec{\chi& c\\  0& 1}$, where $\sigma (q_n)=\epsilon_n^{c_n(\sigma)}q_n$ (cf. \cite[\S~II.4]{Berger}).
 
  The inclusion $\{q_n^s:\ s,n\in\NN	\}\subset \calO_{\bar K}$ gives a pre-log-structure $\bar{N}$ on $\calO_{\bar K}$, and we define $ \Omega^{\log}=\Omega_{(\calO_{\bar K},\bar{N})/(\calO_{K},\text{triv})}$  (see \S~\ref{logstr}). 
  
For any $\ZZ_p$-module $M$ we denote  $M(1)^{\log}=M\otimes_{\ZZ_p}\ZZ_{p}(1)^{\log}$ (and similarly without ${}^{\log}$).  

We aim to prove the following result (cf. {\cite[Théorème 1']{Fon:81a}}).
\begin{theorem}
	There is a natural surjection 
	\[
		\dlog:\bar{K}(1)^{\log}\rightarrow \Omega^{\log}
	\]
	compatible with $\dlog_F$. Moreover it induces an isomorphism 
	\[
		\calO_{\CC_p}(1)^{\log}\cong T_p(\Omega^{\log}) \ .
	\]
\end{theorem}
	In particular $V_p(\Omega^{\log})$ can be seen as a submodule of $\BB_{\rm st}$ sitting in the following exact sequence
\[
	0\to \CC_p\cdot t\to V_p(\Omega^{\log})\to \CC_p\cdot u\to 0
\]
where $t=\log([\underline{\epsilon}^\flat])$ and $u=\log([\underline{q}^\flat])$ (we write $\underline{a}^\flat$ to denote an element of $\CC_p^\flat$ given by a fixed compatible system $\underline{a}$ of $p$-powers roots of $a$).

The result of Fontaine is used to compute the $p$-adic periods of abelian varieties with good reduction. We believe that the logarithmic version  can be used to address the case of semistable abelian varieties in a direct way (see remark~\ref{rmk:pairing}).

\section{Logarithmic differentials} We recall some definitions from \cite{Kat:87}. A logarithmic ring $(A,M,\alpha)$ is the data of: a ring $A$; a monoid $M$; a morphism $\alpha:M\to A$ of monoids\footnote{w.r.t. the multiplication on $A$.}  inducing an isomorphism $\alpha^{-1}(A^\times)\cong A^\times$. If $M$ is a monoid, we denote by $M^+$ its group completion. Given a morphism of logarithmic rings $(A,M,\alpha)\to (B,N,\beta)$ we define the module of log differentials $\Omega_{(B,N)/(A,N)}$ to be the quotient of the module
\[
	\Omega_{B/A}\oplus \left(B\otimes \Coker(M^+\to N^+) \right)
\]
by the submodule generated by elements of the form $(d\beta(n),0)-(0,\beta(n)\otimes n)$ for $n\in N $.
There are natural maps $d:B\to \Omega_{(B,N)/(A,N)}$ (the usual differential) and $\dlog_N:N\to \Omega_{(B,N)/(A,N)}$ such that  $\beta(n)\dlog_N(n)=d\beta(n)$, for all $n\in N$.

\subsection{The log structure on $\calO_{\bar K}$} \label{logstr}
Fix $q\in p\calO_K$, we also fix a system $\underline{q}=(q_n)_n$ of compatible $p$-power roots of $q$: i.e. $q_0=q$ and $q_{n+1}^p=q_n$ for all $n$.
We define the following monoid 
\[
	\bar{N}:= \underline{q}^{\NN[1/p]}\ ,\text{where}\ \underline{q}^{m/p^s}=  q_s^m\ .
\]
Then we have the pre-log-structure   $\bar{N}\to \calO_{\bar K}$.
On $\calO_K$ we may consider pre-log-structure $N=1$.
\subsection{The map $\dlog$}
By definition of $\Omega^{\log}$ there is a natural map $\dlog:\mu_{p^\infty}^{\log}\to \Omega^{\log}$ sending $\epsilon_r^i\cdot q_s^j$ to the image of
\[
	\left(\frac{d\epsilon_r^i}{\epsilon_r^i},1\otimes q_s^j\right)\ .
\]
Clearly we have $\dlog(\epsilon_r\cdot q_s)=\dlog_F(\epsilon_r)+\dlog_{\bar N}(q_s)$.

We denote by the same symbol the $\calO_{\bar K}$-linearisation. Notice that we can identify $\calO_{\bar K}\otimes_{\ZZ_p}\mu_{p^\infty}^{\log} \cong (\bar{K}/ \calO_{\bar{K}})\otimes_{\ZZ_p}\ZZ_{p}(1)^{\log}$ so the previous map can be naturally lifted to
\[
	\dlog:\bar{K}(1)^{\log}\rightarrow \Omega^{\log} \ .
\]
The previous map is compatible with the $\dlog_F$ defined by Fontaine: namely there is a morphism of exact sequences
	$$
	\xymatrix{
	0\ar[r] &\frak{a}(1)\ar[r]\ar[d] &  \bar{K}(1)\ar[r]^{\dlog_F}\ar[d] & \Omega \ar[d] \ar[r]&0\\
	0\ar[r] & Z\ar[r] & \bar{K}(1)^{\log} \ar[r]^{\dlog} & \Omega^{\log}  \ar[r] & 0}
	$$
where $\frak{a} = \{a\in \bar{K} \ :\  v(a)\ge -v(D_{K/K_0})-(p-1)^{-1} \}$ (\cite[Théorème 1 (ii)]{Fon:81a}) and $Z=\ker(\dlog)$.

Now we are ready to complete the proof of the theorem.
\begin{proof}[Proof of the Theorem]
  Consider the following commutative diagram 
  \begin{equation}\label{diag}
  	\xymatrix{
	 {\frak a}\tilde{\epsilon}\ar[r]\ar[d] & \bar{K}\tilde{\epsilon}\ar[r]\ar[d]& \Omega\ar[d]\\
	 Z\ar[d]\ar[r] & \bar K
(1)^{\log}\ar[r]\ar[d]& \Omega^{\log}\ar[d]\\
	  Z/\frak{a}\tilde{\epsilon} \ar[r]&\bar{K}\tilde{q}\ar[r]&\Omega^{\log}/\Omega
	 }
  \end{equation}
	The first two rows (and the last two columns) of \eqref{diag} are  short exact sequences by construction. Hence by the snake lemma the last row and the first column are short exact sequences too.
	
	We now prove that   $(Z/\frak{a}\tilde{\epsilon})\cong\frak{b}\tilde{q}$, where $\frak{b}=\calO_{\bar{K}}$ if $v(q)\ge p$, otherwise
	\[
		\frak{b}=\{b\in \bar{K}\ : \ v(b)\ge \frac{v(q)}{p}-1\} \qquad  ( v(q)< p) \ .
	\]
	Let $b=x/p^n$ with $x\in \calO_{\bar{K}}$ and $n\ge 1$. Then $\dlog(b\otimes \underline{q})$ is the image of $x\otimes q_n $: it maps to $0$ $\Omega^{\log}/\Omega$ if and only if $q_n$ divides $x$, i.e.
	\[
	v(x)\ge v(q_n)=\frac{q}{p^n}\ \Leftrightarrow v(b)\ge v(q_n)=\frac{v(q)}{p^n}-n \ .
	\]
	The maximum of $v(q)p^{-n}-n$ is attained at $n=1$. Since $\dlog(\calO_{\bar{K}} \otimes \underline{q})$ is trivial modulo $\Omega$ the claim is proved.
	
	Now we can deduce the short exact sequence 
	\[
	0\to  (\bar{K}/{\frak a}) \tilde{\epsilon} \to \Omega^{\log} \to   (\bar{K}/\hat{\frak b}) \tilde{q}\to 0 \ 
	\]
	inducing a short exact sequence of Tate modules. Since we have an isomorphism of $\ZZ_p$-modules (for $\frak c =\frak a, \frak b$)
	$$	
		(\bar{K}/{\frak c})[p^n]\cong(p^{-n}{\frak c}/{\frak c})\cong ({\frak c}/p^n{\frak c})\cong (\calO_{\bar K}/p^n)
	$$
	we can conclude the proof of the theorem. Indeed it is possible to check the Galois action on generators and compare with  $\calO_{\CC_p}(1)^{\log}$.
\end{proof}
\begin{remark}\label{rmk:pairing}
 It is now easy to define a  Galois equivariant pairing
 \begin{equation}\label{eq:pair}
	 \langle \cdot,\cdot \rangle:T_pE_q\times \Fil^1H^1_{\dR}(E_q/K)\rightarrow V_p\Omega^{\log}
 \end{equation}

	analogous to the Fontaine paring for abelian varieties with good reduction. We will discuss the details of this construction in a forthcoming work: we will develop  the $p$-adic Hodge theory for $1$-motives,  and construct a perfect pairing between the Tate module and the full de Rham cohomology of a $1$-motive over $K$. The theory in the good reduction (or crystalline) case is already explained in \cite[\S~3.1]{Maz}, we will treat also the (potentially) semi-stable case. 
	
	For the sake of the reader we give a sketch of the construction of \eqref{eq:pair} in the following paragraphs.

	Let $M_q=[\ZZ\xrightarrow{1\mapsto q} \GG_m]$ be the  strict 1-motive over $K$ such that $T_p(M_q) = \ZZ_p(1)^{\log}$ \cite[\S~4.2]{Ray} (See also \cite[\S~1.3.2]{BCC}).
	There is a short exact sequence
	\[
	0\to T_p\Omega^{\log}\to A_2^{\log}\xrightarrow{\theta} \calO_{\CC}\langle X^{\pm 1/p^\infty}\rangle\to 0
	\]
	where $A_2^{\log}$ is simply the $p$-adic completion of $(B_{\dR}^+/\Fil^2)[X^{\pm 1/p^\infty}]$; $\theta$ is the usual map on $B_{\dR}^+$ and $\theta(X)=0$.

	This induces a short exact sequence (of complexes)
	\[
	0\to \Lie(\GG_m)\otimes T_p\Omega^{\log}\to M_q(A_2^{\log})\to M_q(\calO_{\CC}\langle X^{\pm 1/p^\infty}\rangle)\to 0 \ .
	\]
	By taking the cones of the multiplication by $p^n$ and the boundary map of the associated long exact sequence we have
	\[
	\rho_n:  M_q(\calO_{\CC})[p^n]=(\bar{K}^*/\langle q\rangle)[p^n]\rightarrow \Lie(\GG_m)\otimes (T_p\Omega^{\log}/p^n) \ ,
	\]
	since $M_q(\calO_{\CC}\langle X^{\pm 1/p^\infty}\rangle)[p^n]=M(\calO_{\CC})[p^n]$.
	Then taking the limit over $n$ induces a pairing
	\[
	\langle \cdot,\cdot \rangle:\ZZ_p(1)^{\log}\times \Fil^1H^1_{\dR}(\GG_m/\calO_K)\rightarrow T_p\Omega^{\log}
\]
such that $\langle (x_n)_n,dT/T \rangle=(\dlog(x_n))_n$. To get the claim it is sufficient to invert $p$ and use the isomorphisms given by the rigid uniformisation:  $\Fil^1H^1_{\dR}(\GG_m/K)\cong \Fil^1H^1_{\dR}(E_q/K)$ and  $\ZZ_p(1)^{\log}\cong T_pE_q$.
\end{remark}

\end{document}